\documentclass[12pt,twoside]{article}    
\usepackage{amssymb,
		amsmath,
		dsfont,
		latexsym,
		enumerate,
		xspace,
		graphicx,
		rotating,
		mathrsfs,
		color,
		psfrag,
		ifpdf,
	}
\usepackage[utf8]{inputenc}
\usepackage{theorem}

\makeatletter
\@addtoreset{equation}{section}
\makeatother

\newtheorem{theorem}{Theorem}[section]

\newtheorem{definition}[theorem]{Definition}

\newtheorem{lemma}[theorem]{Lemma}
\newtheorem{Remarks}[theorem]{Remarks}

\newcommand{\BOX}{\mbox{{\ensuremath{\Box}}\hspace{-0.5mm}}}

\newcommand{\conv}{{\mathop{\mathrm{conv}}}}

\newcommand{\dist}{{\mathop{\mathrm{dist}}}}
\newcommand{\proj}{{\mathop{\mathrm{proj}}}}

\newcommand{\SekEins}[1]{\relax}

\newcommand{\N}{\mathds N}

\newcommand{\R}{\mathds R}

\let\EINS\1
\newcommand{\1}{\mathds 1}

\newenvironment{proof}{{\vskip\baselineskip\noindent\textbf{Proof:}}}%
{\origqed\vskip\baselineskip\gdef\origqed{\hspace*{.1pt}\hspace*{\fill}\BOX}}
\newenvironment{proofx}[1]{{\vskip\baselineskip\noindent\textbf{Proof #1:}}}%
{\origqed\vskip\baselineskip\gdef\origqed{\hspace*{.1pt}\hspace*{\fill}\BOX}}
\renewcommand{\BOX}{\ensuremath\Box}
\makeatletter
\def\@tagformdelimstart{(}%
\def\@tagformdelimend{)}%
\def\@tagformdel{%
 \gdef\@tagformdelimstart{}%
 \gdef\@tagformdelimend{}%
}
\def\@tagformset{%
 \gdef\@tagformdelimstart{(}%
 \gdef\@tagformdelimend{)}%
}
\def\tagform@#1{%
  \maketag@@@{\@tagformdelimstart\ignorespaces#1\unskip%
  \@@italiccorr\@tagformdelimend}\@tagformset}
\def\origqed{\hspace*{.1pt}\hspace*{\fill}\BOX}
\def\qed{\ifmmode%
 \@tagformdel%
 \tag{\BOX}%
 \else%
 \hspace*{.1pt}\hspace*{\fill}\BOX%
 \fi%
 \gdef\origqed{}}
\makeatother
\def\today{\number\day\space \ifcase\month \or Jan \or Feb \or Mar \or Apr \or
Mai \or Jun \or Jul \or Aug \or Sep \or Oct \or Nov \or Dec \fi \number\year}

\newcount\hours
\newcount\minutes

\def\SetTime{\hours=\time
    \global\divide\hours by 60
    \minutes=\hours
    \multiply\minutes by 60
    \advance\minutes by-\time
    \global\multiply\minutes by-1 }
\SetTime

\def\now{\number\hours:\ifnum\minutes<10 0\fi\number\minutes}

\makeatother

\title{\bf \Large Blackwell Prediction for Categorical Data}
\author{H. R. Lerche \\ University of Freiburg i. Br.}
\date{ }

\begin{document}
\maketitle
\thispagestyle{empty}


\begin{abstract}
We study the problem of sequential prediction of categorical data and discuss a generalisation of Blackwell's
algorithm on 0-1 data. The arguments are based on Blackwell's approachability results given in \cite{Blackwell}.
They use mainly linear algebra.
\end{abstract}

\footnotetext{\textit{Date:} \date{\today}}

\section{Introduction and Background}%
\label{sec:1}

Let us consider the problem of sequential prediction of categorical data. Let $D=\{0,1, \ldots, d-1 \}$ denote the set of possible outcomes with $d \geq 2$. Let $x_1, x_2, \ldots$ be an infinite sequence with values in $D$. Let $Y_1, Y_2, \ldots$ denote the sequence of predictions.
This is a random sequence with values in $D$. $Y_{n+1}$ predicts $x_{n+1}$ and may depend on the first $n$ outcomes $x_1, x_2, \ldots, x_n, Y_1,Y_2, \ldots, Y_n$ and some additional random mechanism. Our goal ist to construct a sequential prediction procedure which works well for all sequences $(x_i)_{i \in \N}$ in an asymptotic sense. %
We intend to generalize Blackwell's prediction procedure for two categories. The algorithm of Blackwell can be described as follows using Figure \ref{fig:1} below. Let $x_1, x_2,\dots$ be an infinite 0-1 sequence. Let
$\overline{x}_n =\frac1n\sum_{k=1}^n x_k$ be the relative frequency of the ``ones'' and $\overline{\gamma}_n=\frac1n\sum_{k=1}^n \1_{\{Y_k = x_k\}}$ the relative frequency of correct guesses. Let $\mu_n=(\overline x_n,\overline{\gamma}_n)\in[0,1]^2$ and $\mathcal S = \{(x,y)\in[0,1]^2\mid y\ge \max(x,1-x)\}$.
%
\begin{center}
\begin{picture}(0,0)%
\includegraphics{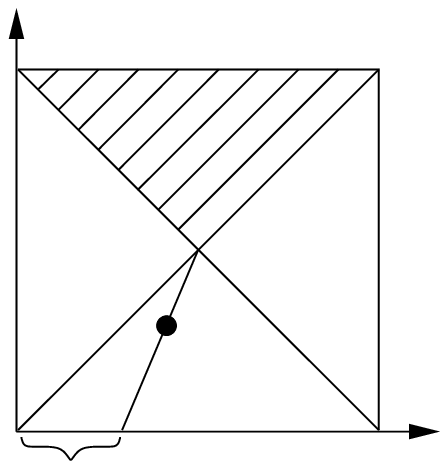}%
\end{picture}%
\setlength{\unitlength}{2072sp}%
\begingroup\makeatletter\ifx\SetFigFontNFSS\undefined%
\gdef\SetFigFontNFSS#1#2#3#4#5{%
  \reset@font\fontsize{#1}{#2pt}%
  \fontfamily{#3}\fontseries{#4}\fontshape{#5}%
  \selectfont}%
\fi\endgroup%
\begin{picture}(4399,4556)(1744,-5585)
\put(2670,-3312){\makebox(0,0)[lb]{\smash{{\SetFigFontNFSS{12}{14.4}{\familydefault}{\mddefault}{\updefault}{\color[rgb]{0,0,0}$D_1$}%
}}}}
\put(3870,-4595){\makebox(0,0)[lb]{\smash{{\SetFigFontNFSS{12}{14.4}{\familydefault}{\mddefault}{\updefault}{\color[rgb]{0,0,0}$D_3$}%
}}}}
\put(4787,-3309){\makebox(0,0)[lb]{\smash{{\SetFigFontNFSS{12}{14.4}{\familydefault}{\mddefault}{\updefault}{\color[rgb]{0,0,0}$D_2$}%
}}}}
\put(3757,-4022){\makebox(0,0)[lb]{\smash{{\SetFigFontNFSS{12}{14.4}{\familydefault}{\mddefault}{\updefault}{\color[rgb]{0,0,0}$\mu_n$}%
}}}}
\put(5453,-5329){\makebox(0,0)[lb]{\smash{{\SetFigFontNFSS{12}{14.4}{\familydefault}{\mddefault}{\updefault}{\color[rgb]{0,0,0}$\overline{x}_n$}%
}}}}
\put(2646,-5462){\makebox(0,0)[lb]{\smash{{\SetFigFontNFSS{12}{14.4}{\familydefault}{\mddefault}{\updefault}{\color[rgb]{0,0,0}$w_n$}%
}}}}
\put(1759,-1768){\makebox(0,0)[lb]{\smash{{\SetFigFontNFSS{12}{14.4}{\familydefault}{\mddefault}{\updefault}{\color[rgb]{0,0,0}$\overline{\gamma}_n$}%
}}}}
\put(3871,-2221){\makebox(0,0)[lb]{\smash{{\SetFigFontNFSS{17}{20.4}{\familydefault}{\mddefault}{\updefault}{\color[rgb]{0,0,0}$\mathcal{S}$}%
}}}}
\end{picture}%
\label{fig:1} ~\\
Figure 1
\end{center}
%

In Fig. \ref{fig:1}, let $D_1$, $D_2$ and $D_3$ be the left, right, and bottom triangles, respectively, in the unit square so that $D_1 = \{(x,y) \in [0,1]^2 \mid x \leq y \leq 1-x \}$ etc. When $\mu_n \in D_3$, draw the line through the points $\mu_n$ and $(\frac12 , \frac12)$ and let $(w_n , 0)$ be the point where this line crosses the horizontal axis. The Blackwell algorithm chooses its prediction $Y_{n+1}$ on the basis of $\mu_n$ according to the (conditional) probabilities
\begin{equation*}
 P(Y_{n+1} =1) = \left\{ \begin{array}{ll}
          0 & \mbox{if } \mu_n \in D_1 \\
          1 & \mbox{if } \mu_n \in D_2 \\
          w_n & \mbox{if } \mu_n \in D_3.
         \end{array} \right.
\end{equation*}
When $\mu_n$ is in the interior of $\mathcal{S}$, $Y_{n+1}$ can be chosen arbitrarily. Let $Y_1$ = 0. It then holds that for the Blackwell algorithm applied to any 0-1 sequence $x_1, x_2, \ldots$ the sequence $(\mu_n ; \ n \geq 1)$ converges almost surely to $\mathcal{S}$, i.e. $\dist (\mu_n, \mathcal{S}) \to 0$ as $n \to \infty$ almost surely. Here $\dist(\cdot, \cdot)$ denotes the Euclidean distance from $\mu_n$ to $\mathcal{S}$.

As Blackwell once pointed out this is a direct consequence of his Theorem~1 in \cite{Blackwell} when one chooses the payoff matrix as
$$
\begin{pmatrix}
(0,1) & (1,0) \\ (0,0) & (1,1)
\end{pmatrix} .
$$
For a quick almost sure argument see \cite{RobbinsSiegmund}.
Blackwell also raised the question whether his Theorem~1 of \cite{Blackwell} applies to sequential prediction when there are more than two categories. We shall study this question and finally answer it affirmative.

We construct a Blackwell type prediction procedure for $d>2$ categories by choosing the state space and the randomisation rules in a certain way. This procedure then has similar properties as Blackwell's original one. It also has the feature that the $d$-category procedure reduces to the $(d-1)$ category procedure if one category is not observed.

The structure of this paper is as follows. In Section \ref{sec:2} we introduce the appropriate state space and define the randomisation rule. In Section \ref{sec:3} we state the convergence result and prove it. For that we shall apply a simplified version of Blackwell's Theorem 1 of \cite{Blackwell}, which we also state in Section \ref{sec:3}. 

This paper is a continuation of \cite{LercheSakar}, where the case $d=3$ was discussed, and of the diploma thesis of R. Sandvoss \cite{Sandvoss}.

We shall use the following notation: Latin letters for points, vectors, and indices, greek letters for scalars. We denote components of vectors or points by superindices like $v=(v^{(0)}, \dots, v^{(d-1)})\in\R^d$. $e_0=(1,0,\dots,0), \dots, e_{d-1}=(0,\dots,0,1)$ denote the $d$-dimensional unit points and $\1_d=(1,\dots,1)$. The affine subspace of $\R^d$ generated by the points $a_0,\dots, a_n\in\R^d$ is given by
\[
A(\{a_0,\dots,a_n\}) \kern-1pt:=\kern-1pt \bigg\{\kern-1pt a\in\R^d \,{\Big|}\,
 a = \kern-2pt\sum_{i=0}^n \kern-.5pt\lambda_i a_i,  \sum_{i=0}^n \kern-.5pt\lambda_i=1, \lambda_i\in\R, a_i\in\R^d,i=0,\dots,n\bigg\}\kern-.5pt.
\]
The convex hull of $a_i,\dots,a_n\in\R^d$ is given by
\[
\begin{array}{ll}
\lefteqn{\conv(\{a_0,\dots,a_n\})} \\
\quad & = \bigg\{a\in\R^d \,{\Big|}\,
a = \sum_{i=0}^n \lambda_i a_i, \sum_{i=0}^n\lambda_i=1, \lambda_i\in[0,1], a_i\in\R^d,i=0,\dots,n\bigg\}.
\end{array}
\]
%
The Euclidean scalar product on $\R^d$ is given by $\langle\cdot,\cdot\rangle$, the Euclidean distance by $\dist(\cdot,\cdot)$.


\section{The Construction of the \boldmath$d$-Dimensional Prediction Procedure}%
\label{sec:2}
\subsection{The Structure of the Prediction Prism}%
\label{subsec:2.1}

For $n\in\N$, $x_1, x_2,\dots,x_n\in D$ let $Y_1,Y_2,\dots,Y_n\in D$ denote the corresponding predictions. Let $\overline x_n=(\overline x_n^{(0)},\dots,\overline x_n^{(d-1)})$ with $\overline x_n^{(l)}=\frac1n\sum_{k=1}^n \1_{\{x_i=l\}}$, $l\in D$, denote the vector of the relative frequencies of the $n$ outcomes and $\overline\gamma_n =\frac1n\sum_{k=1}^n \1_{\{Y_k=x_k\}}$ the relative frequency of correct predictions.

Let
$$
\mbox{\boldmath $\Sigma$\unboldmath$_{d-1}$} =\left\{ (q_0,\dots,q_{d-1}) \mid q_l\ge 0, \ \mbox{$\sum\limits_{l=0}^{d-1}q_l=1$} \right\} $$
denote the unit simple in $\R^d$ and
$$W_d=\mbox{\boldmath $\Sigma$\unboldmath$_{d-1}$} \times [0,1] = \left\{ (q,\gamma) \mid q\in \Sigma_{d-1}, \ 0\le \gamma\le 1 \right\}.$$

Since $\sum_{l=0}^{d-1} x_n^{(l)}=1$, we have
$\overline x_n \in \mbox{\boldmath $\Sigma$\unboldmath$_{d-1}$} $
and
$(\overline x_n,\overline \gamma_n)\in W_d$.
Let $\mathscr{S}_d=\{(q,\gamma)\in W_d \mid \gamma\ge \max_l q^{(l)}\}$. We are interested in prediction procedures for which $\mu_n:=(\overline x_n,\overline \gamma_n)$ converges to $\mathscr{S}_d$ for every sequence $x_1,x_2,\dots$ This means that the Euclidean distance $\dist(\mu_n,\mathscr{S}_d)\to 0$ as $n\to\infty$.

Unfortunately Blackwell's Theorem~1 of \cite{Blackwell} cannot be applied directly. The reader may take a look at Theorem~\ref{theo:2} below which is a simplified version of Blackwell's result. The condition (C) there does not hold in general for $W_d$ and $\mathcal{S}_d$. (To see this, let $d=3$, $s=(\frac13,\frac13,\frac13,\frac13)$, $\mu_n=(\frac14,\frac14,\frac12,0)$. Then $p(\mu_n)=\mu_n$, and $s-\mu_n$ is not perpendicular to $\mathcal R(p(\mu_n))$.)

The difficulties vanish when one modifies the state space in the right way. Let
$V_d=\{q+\gamma\1_d\mid (q,\gamma)\in W_d\}$ with $\1_d=(1,\dots,1)$.
Then $v_n:=\overline x_n+\overline \gamma_n\1_d\in V_d$ for all $n$.
The convergence of $\mu_n$ to $\mathscr{S}_d$ corresponds to that of $v_n$ to
\boldmath$\mathcal{S}$\unboldmath$_d$ where
\boldmath$\mathcal{S}$\unboldmath$_d=\{q+\gamma\1_d\in V_d \mid \gamma\ge \max_{l} q^{(l)}\}$. This follows from the fact that $\Psi:W_d\to V_d$ with $\Psi((q,\gamma))=q+\gamma\1_d$ is an isometric bijection of $W_d$ on $V_d$. We note that for $z,z'\in W_d$ it holds that
\[
\dist(\Psi(z),\Psi(z'))^2 = \sum_{i=0}^{d-1} (z_i-z'_i)^2+d\cdot (z_d-z'_d)^2.
\]

To construct the appropriate randomisation regions let us ``cut'' the prism $V_d$ by certain hyperplanes. (This corresponds to splitting the unit square by the diagonals in the case of two categories.)

Let $e_0=(1,0,0,\dots,0),\dots, e_{d-1}=(0,0,\dots,0,1)$ denote the $d$-dimensional unit points. Let
$E_l=A(\{e_0,\dots,e_{l-1},e_l+\1_d,e_{l+1},\dots, e_{d-1}\})$, $l=0,\dots,d-1$, denote the hyperplanes which contain one vertex of the ``upper side'' of the prism $e_l+\1_d$ and $(d-1)$ vertices $e_k\not= e_l$ of
\boldmath$\mathcal{S}$\unboldmath$_{d-1}$. The $d$ hyperplanes $E_l$ cut the prism $V_d$ in $2^d$ pieces, and all contain the point $s=(\frac2d,\frac2d,\dots,\frac2d)$. In this point $s$ the planes $E_l$ are all perpendicular to each others.

This can easily be seen since their corresponding normal vectors are given by
$n_l=-e_l+\frac2d \1_d$. This leads to the following characterization of lying ``above'' $E_i$:
\[
v \text{ lies above } E_i\Leftrightarrow \langle v-n_i,n_i\rangle < 0.
\]
In the same way one defines lying below and in $E_i$.

Now we can describe \boldmath$\mathcal{S}$\unboldmath$_d$ in two different ways:
\begin{eqnarray*}
\mathcal S_d &=& \{q+\gamma\1_d\in V_d\mid \langle q-n_l,n_l\rangle \ge 0 \text{ for } l=0,\dots,d-1\}\\
&=& \{q+\gamma\1_d\in V_d\mid \gamma \ge \max(q^{(0)},\dots,q^{(d-1)})\}.
\end{eqnarray*}

For the case $d=3$ the sets $V_d$ and $\mathcal S_d$ are shown in the following figures.

\begin{minipage}{.48\textwidth}%
\begin{center}
\begin{picture}(0,0)%
\includegraphics{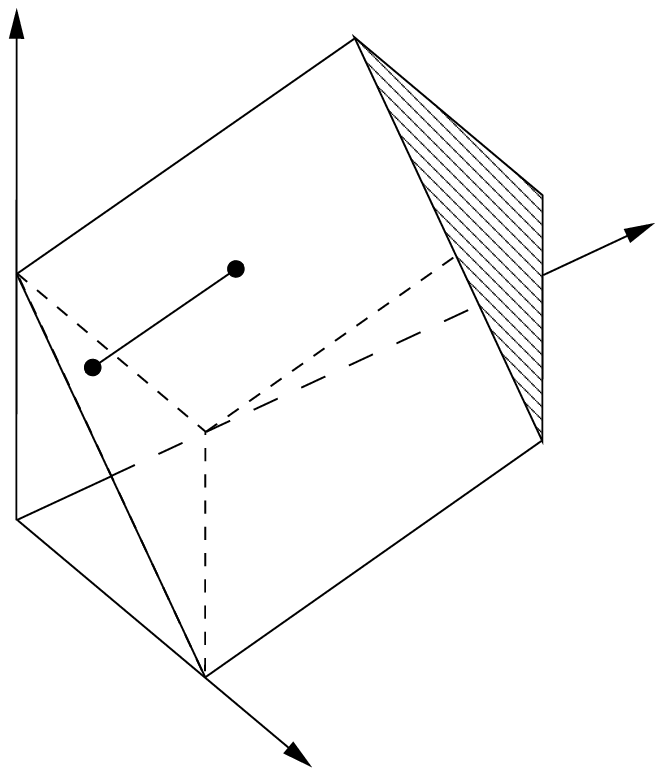}%
\end{picture}%
\setlength{\unitlength}{2072sp}%
\begingroup\makeatletter\ifx\SetFigFontNFSS\undefined%
\gdef\SetFigFontNFSS#1#2#3#4#5{%
  \reset@font\fontsize{#1}{#2pt}%
  \fontfamily{#3}\fontseries{#4}\fontshape{#5}%
  \selectfont}%
\fi\endgroup%
\begin{picture}(6681,7259)(-374,-6668)
\put(5581,-1321){\makebox(0,0)[lb]{\smash{{\SetFigFontNFSS{12}{14.4}{\familydefault}{\mddefault}{\updefault}{\color[rgb]{0,0,0}$x^{(2)}$}%
}}}}
\put(3421,-5281){\makebox(0,0)[lb]{\smash{{\SetFigFontNFSS{12}{14.4}{\familydefault}{\mddefault}{\updefault}{\color[rgb]{0,0,0}$V_{3}$}%
}}}}
\put(1981,-6541){\makebox(0,0)[lb]{\smash{{\SetFigFontNFSS{12}{14.4}{\familydefault}{\mddefault}{\updefault}{\color[rgb]{0,0,0}$x^{(1)}$}%
}}}}
\put(1621,-4201){\makebox(0,0)[lb]{\smash{{\SetFigFontNFSS{12}{14.4}{\familydefault}{\mddefault}{\updefault}{\color[rgb]{0,0,0}$\Sigma_2$}%
}}}}
\put(2611,-2221){\makebox(0,0)[lb]{\smash{{\SetFigFontNFSS{12}{14.4}{\familydefault}{\mddefault}{\updefault}{\color[rgb]{0,0,0}$v_n$}%
}}}}
\put(1261,-3211){\makebox(0,0)[lb]{\smash{{\SetFigFontNFSS{12}{14.4}{\familydefault}{\mddefault}{\updefault}{\color[rgb]{0,0,0}$\overline{x}_n$}%
}}}}
\put(-359,119){\makebox(0,0)[lb]{\smash{{\SetFigFontNFSS{12}{14.4}{\familydefault}{\mddefault}{\updefault}{\color[rgb]{0,0,0}$x^{(0)}$}%
}}}}
\end{picture}%

\label{fig:2}~ \\[-1ex]
Figure 2
\end{center}
\end{minipage}%
\hfill
\begin{minipage}{.48\textwidth}%
\begin{center}
~\\[5mm]
\begin{picture}(0,0)%
\includegraphics{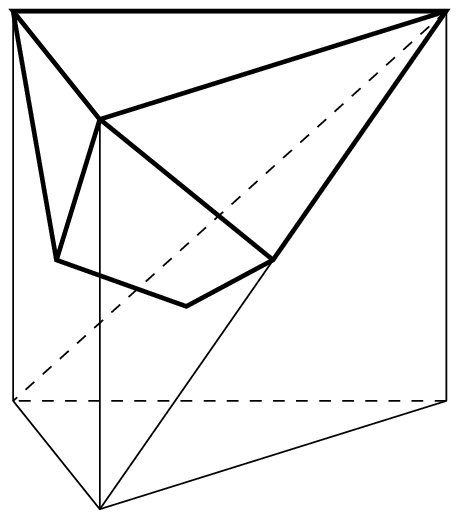}%
\end{picture}%
\setlength{\unitlength}{1823sp}%
\begingroup\makeatletter\ifx\SetFigFontNFSS\undefined%
\gdef\SetFigFontNFSS#1#2#3#4#5{%
  \reset@font\fontsize{#1}{#2pt}%
  \fontfamily{#3}\fontseries{#4}\fontshape{#5}%
  \selectfont}%
\fi\endgroup%
\begin{picture}(5160,5700)(976,-7501)
\put(3016,-5326){\makebox(0,0)[lb]{\smash{{\SetFigFontNFSS{11}{13.2}{\rmdefault}{\mddefault}{\updefault}{\color[rgb]{0,0,0}$s$}%
}}}}
\put(991,-6046){\makebox(0,0)[lb]{\smash{{\SetFigFontNFSS{11}{13.2}{\familydefault}{\mddefault}{\updefault}{\color[rgb]{0,0,0}2}%
}}}}
\put(5761,-6361){\makebox(0,0)[lb]{\smash{{\SetFigFontNFSS{11}{13.2}{\familydefault}{\mddefault}{\updefault}{\color[rgb]{0,0,0}1}%
}}}}
\put(2161,-7486){\makebox(0,0)[lb]{\smash{{\SetFigFontNFSS{11}{13.2}{\familydefault}{\mddefault}{\updefault}{\color[rgb]{0,0,0}0}%
}}}}
\put(6121,-2176){\makebox(0,0)[lb]{\smash{{\SetFigFontNFSS{12}{14.4}{\familydefault}{\mddefault}{\updefault}{\color[rgb]{0,0,0}$\mathcal{S}$}%
}}}}
\end{picture}%

\label{fig:3} ~\\[15mm]
Figure 3
\end{center}
\end{minipage}

\bigskip

\subsection{The Randomisation Rule}%
\label{subsec:2.2}

For $v_n=\overline x_n+\overline\gamma_n\1_d$ we will define a $d$-dimensional random vector $p(v_n)\in \mbox{\boldmath $\Sigma$\unboldmath$_{d-1}$}$. It plays the same role as $w_n$ does in the 0-1 case. With it we define $Y_{n+1}:$
$$P(\{Y_{n+1}=k\})=p^{(k)}(v_n) \mbox{ for } k\in D.$$
\begin{definition}\label{def:1}
Let $v_n\in V_d$, $n\in \N$ and let $(i_0,\dots,i_{d-1})$ be a permutation of $(0,\dots,d-1)$ such that it holds:
\begin{eqnarray*}
&& \langle v_n-n_l,n_l\rangle \le 0\quad \text{for } l=i_0,\dots,i_j\\
\text{and}&& \langle v_n-n_l,n_l\rangle > 0 \quad \text{for } l=i_{j+1},\dots,i_{d-1}.
\end{eqnarray*}
\end{definition}

\noindent%
\textbf{Case 1:} %
Let $v_n\in V_d\setminus \mathcal S_d$.

Let $A_1=A(\{\frac2d\1_d,e_{i_{j+1}}, \dots, e_{i_{d-1}}, v_n\})$ be the affine space of $\R^d$ generated by the points in the waved brackets.
Let $A_2=A(\{ e_{i_0},\dots, e_{i_j}\})$ denote the corresponding affine space. The intersection $A_1\cap A_2$ contains exactly one point of \boldmath $\Sigma$\unboldmath$_{d-1}$, we call it $p(v_n)$.

\noindent%
\textbf{Case 2:} %
Let $v_n \in \partial \mathcal S_d$. Let $\nu = \#\{E_k\mid v_n\in E_k \text{ for } k=0,\dots,d-1\}$.

Then
\[
p^{(k)}_{(v_n)} = \begin{cases}
         1/\nu &\text{ for } v_n\in E_k\\
		 0 & \text{ for } v_n\not\in E_k
         \end{cases}
\]
for $k=0,1,\dots,d-1$.

The prediction procedure just defined is called ``Generalized Blackwell algorithm''.
\begin{Remarks}
\begin{enumerate}[1)]
\item %
The case $v_n\in\mathcal S_d\setminus \partial \mathcal S_d$ does not occur by the construction of the rule.
\item %
$A_2=\emptyset$ cannot occur, since then there exists at least one $k\in D$ with
\hbox{$\langle v_n-n_k,n_k\rangle\le 0$}.
\item %
We note that $A_1\cap A_2$ contains always just one point of \boldmath $\Sigma$\unboldmath$_{d-1}$.
\item %
For $j=d-1$ one obtains %
$A_1= A(\{\frac2d \1_d,v_n\})$, %
$A_2= A(\{e_{i_0},\dots, e_{i_{d-1}}\})$ and $p(v_n)$ is the projection along the line, defined by $\frac2d\1_d$ and $v_n$ ``down'' to $\Sigma_{d-1}$.
\item %
For $d=3$ the following figure shows the randomisation in a ``lower'' side piece of the prism. Here planes lie above $\mu_n$ and one below.
\begin{center}
\begin{picture}(0,0)%
\includegraphics{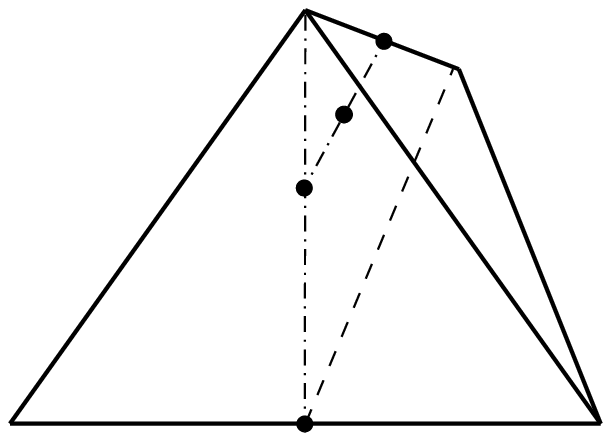}%
\end{picture}%
\setlength{\unitlength}{2486sp}%
\begingroup\makeatletter\ifx\SetFigFontNFSS\undefined%
\gdef\SetFigFontNFSS#1#2#3#4#5{%
  \reset@font\fontsize{#1}{#2pt}%
  \fontfamily{#3}\fontseries{#4}\fontshape{#5}%
  \selectfont}%
\fi\endgroup%
\begin{picture}(5043,4037)(-509,257)
\put(2731,4027){\makebox(0,0)[lb]{\smash{{\SetFigFontNFSS{12}{14.4}{\familydefault}{\mddefault}{\updefault}{\color[rgb]{0,0,0}nearest point in $\mathcal{S}$ for $v_n$}%
}}}}
\put(3601,3449){\makebox(0,0)[lb]{\smash{{\SetFigFontNFSS{14}{16.8}{\familydefault}{\mddefault}{\updefault}{\color[rgb]{0,0,0}$s$}%
}}}}
\put(1936,389){\makebox(0,0)[lb]{\smash{{\SetFigFontNFSS{14}{16.8}{\familydefault}{\mddefault}{\updefault}{\color[rgb]{0,0,0}$p(v_n)$}%
}}}}
\put(-494,389){\makebox(0,0)[lb]{\smash{{\SetFigFontNFSS{14}{16.8}{\familydefault}{\mddefault}{\updefault}{\color[rgb]{0,0,0}$(1,0,0)$}%
}}}}
\put(2566,2729){\makebox(0,0)[lb]{\smash{{\SetFigFontNFSS{14}{16.8}{\familydefault}{\mddefault}{\updefault}{\color[rgb]{0,0,0}$v_n$}%
}}}}
\put(4006,389){\makebox(0,0)[lb]{\smash{{\SetFigFontNFSS{14}{16.8}{\familydefault}{\mddefault}{\updefault}{\color[rgb]{0,0,0}$(0,1,0)$}%
}}}}
\end{picture}%

{\rm Figure 4}
\end{center}
\end{enumerate}
\end{Remarks}

\section{The Convergence Result}\label{sec:3}
\subsection{Main Result}%
\label{subsec:3.1}

\begin{theorem}\label{theo:1}
Let $d\ge 2$. Then for the generalized Blackwell algorithm, applied to any infinite sequence $x_1,x_2,\dots$ with values in $D$, it holds that $\dist(v_n, \mathcal{S}_d)\to 0$ with probability one as $n\to\infty$.
\end{theorem}

Now we shall derive Theorem \ref{theo:1} by tracing it back to Blackwell's Theorem~1 of \cite{Blackwell}. This we first state in a simplified version.

\subsection{Blackwell's Minimax Theorem}%
\label{subsec:3.2}

We consider a repeated game of two players with a payoff matrix $M=(m_{ij})$ with $m_{ij}\in\R^d$ and $1\le i\le r$ and $1\le j\le s$. Player I chooses the row, player II the column. Let
\[
\mathcal P=\bigg\{p=(p_1,\dots,p_r) \,{\Big|}\, p_i\ge 0, \sum_{i=1}^r p_i=1\bigg\}
\]
denote the mixed actions of player I and
\[
\mathcal Q=\bigg\{q=(q_1,\dots,q_s) \,{\Big|}\, q_j\ge 0, \sum_{j=1}^s q_j=1\bigg\}
\]
the mixed actions of player II. %
A strategy $f$ in a repeated game for player I is a sequence $f=(f_k; \ k \geq 1)$ with $f_k \in \mathcal{P}$. A strategy $g$ for player II is defined similarly. Two strategies define a sequence of payoffs
$z_k$, $k=1,2, \ldots$
In detail: If in the $k$-th game $i$ and $j$ are choosen according to $f_k$ and $g_k$, the payment to player I is $m_{ij} \in \R^d$. Blackwell discussed in
\cite{Blackwell}
the question: Can player I control $\overline{z}_n = \frac{1}{n} \sum_{k=1}^{n} z_k$ with a certain strategy such that $\overline{z}_n$ approaches a given set $\mathcal{S}$ independently of what player II does?
\begin{definition}\label{def:2}
 A set $\mathcal{S} \subset \R^d$ is approachable for player I if there exists a strategy $f^*$ for which $\dist(\overline{z}_n, \mathcal{S}) \to 0$ with probability one.
\end{definition}

\begin{theorem}[Blackwell] \label{theo:2}
 For $p \in \mathcal{P}$ let
$$
\mathcal{R}(p) = \conv \left( \sum_{i=1}^{r} p_i m_{ij} ; \ j = 1,2, \ldots, s \right).
$$
Let $\mathcal{S}$ denote a closed convex subset of $\R^d$. For every $z \not\in \mathcal{S}$ let $y$ denote the closest point in $\mathcal{S}$ to $z$. We assume:
\begin{itemize}
\item [\rm{(C)}] %
For every $z \not\in \mathcal{S}$ there exists a $p(z) \in \mathcal{P}$ such that the hyperplane through $y$, 
which is perpendicular to the line segment $\overline{zy}$, seperates $z$ from $\mathcal{R}(p(z))$.
\end{itemize}
Then $\mathcal{S}$ is approachable for player I.
\end{theorem}

\subsection{Proof of the Main Result}
\label{subsec:3.3}

To apply Theorem \ref{theo:2} to our case, we choose the vertices of $V_d$ as ``payments'':
\[
m_{ij}=\begin{cases}
    e_i+\1_d & \text{ if } i=j,\\
	e_j & \text{ if } i\not= j.
    \end{cases}
\]
We choose $\mathcal{S}$ as $\mathcal S_d=\{q+\gamma\1_d\in V_d\mid \gamma\ge \max_l q^{(l)}\}$. Then
\begin{eqnarray*}
\mathcal R(p) %
&=& \conv\bigg(\bigg\{\sum_{i=0,i\not=j}^{d-1} p^{(i)} e_j+p^{(j)} (e_j+\1_d)\,\Big|\, j=0,\dots,d-1\bigg\}\bigg) \\
&=& \conv\bigg(\bigg\{\sum_{i=0}^{d-1} p^{(i)} e_j+p^{(j)}\1_d \,\Big|\, j=0,\dots,d-1\bigg\}\bigg)
\\
&=& \conv \big(\big\{ e_j+p^{(j)}\1_d \mid j=0,\dots,d-1 \big\}\big).
\end{eqnarray*}

\noindent
It is left to show that condition (C) is fulfilled.

\noindent
Let $v\in V_d\setminus \mathcal S_d$. We denote by $v_{\proj}$ the closest point in $\mathcal S_d$ to $v$. We will show:

\begin{enumerate}[1)]
\item[\textbf{Fact 1}]$v_{\proj}\in\mathcal R(p(v))$

\item[\textbf{Fact 2}] $v-v_{\proj}$ is perpendicular to $A(\mathcal R(p))$. Here $A(\mathcal R(p))$ means the smallest affine subspace which contains $\mathcal R(p)$.
\end{enumerate}

Both facts together imply condition (C) and finally Theorem \ref{theo:1}.

For the proofs we shall assume that the following situation holds: For $v\in V_d\setminus \mathcal S_d$ it holds
\begin{eqnarray*}
&& \langle v-n_i,n_i\rangle \le 0\quad \text{for } i=0,\dots,j\\
\text{and}&& \langle v-n_i,n_i\rangle > 0 \quad \text{for } i=j+1,\dots,d-1.
\end{eqnarray*}

\begin{proofx}{of Fact 1}
$v$ lies below $E_i$ for $i=0,1,\dots,j$, but $v_{\proj}\in \mathcal S_d$. Thus $v_\proj\in E_0\cap\dots\cap E_j$.
Then
\[
E_0\cap\dots\cap E_j= A\bigg(\bigg\{ e_{j+1}, \dots, e_{d-1}, \frac2d\1_d \bigg\}\bigg).
\]
Thus
\begin{eqnarray*}
v_\proj &\in& A\bigg(\bigg\{e_{j+1},\dots,e_{d-1}\frac2d\1_d \bigg\}\bigg) \cap V_d
\\
&& \subset A\big(\big\{e_i+p^{(i)}(v)\1_d \mid i=0,\dots,d-1 \big\}\big) \cap V_d
= \mathcal R(p(v)).
\end{eqnarray*}
The inclusion follows since $p^{(l)}(v)=0$ for $j+1=l\le d-1$ and $\frac2d\1_d =\frac1d \sum_{i=0}^{d-1} (e_i+p^{(i)}\1_d)$. ~\vspace*{-4mm}
\end{proofx}

Fact 2 will be proven by a sequence of lemmata. At first we generate a new auxiliary point $\tilde v$ which lies in the same plane as $p(v)$.

\begin{lemma}\label{lem:1}
For $v\in V_d\setminus \mathcal S_d$ let $A'=A(\{v,v_\proj\})$ and $A''=A(\{e_{j+1},\dots,e_{d-1},p(v)\})$. Then there exists exactly one point $\tilde v \in A'\cap A''$ and $\tilde v\not\in S_d$.
\end{lemma}

\begin{proof}
Let $A_1=A(\{e_{j+1},\dots,e_{d-1},\frac2d\1_d,v\})$ as in Definition \ref{def:1}. Then according to Definition \ref{def:1} $p(v)\in A_1$ and $v_\proj\in A_1$ by the proof of Fact 1. Then it follows that $\frac2d\1_d\in A'\vee A''$. Here $A'\vee A''$ denotes the smallest affine space, which contains $A'$, $A''$. It holds $A_1=A'\vee A''$.
Since $A'$ and $A''$ are not parallel it follows that $A'\cap A''\not=\emptyset$ and by the dimension formula $\dim(A'\cap A'')= 0$. Hence $A'\cap A''$ contains exactly one point. We call it $\tilde v$. If $\tilde v\in \mathcal S_d$, then $\tilde v\in \mathcal S_d\cap A''$. Then
$\mathcal S_d\cap A(\mbox{\boldmath $\Sigma$\unboldmath}_{d-1})\not=\emptyset$, which is a contradiction to the definitions of $\mathcal S_d$ and \boldmath $\Sigma$\unboldmath$_{d-1}$.
\end{proof}

A direct consequence of Lemma \ref{lem:1} is

\noindent%
\begin{tabbing}\textbf{Fact 3:} \= a) \=
$v_\proj=(\tilde v)_\proj$; \\[1.5ex]
\> b) \>
$v-v_\proj\,\bot\, A(\mathcal R(p))\Leftrightarrow \tilde v-(\tilde v)_\proj \,\bot\, A(\mathcal R(p))$.
\end{tabbing}

We shall use Fact 3 to show Fact 2. At first we calculate $(\tilde v)_\proj$ from $\tilde v$.
For simplification, we write  $\tilde v_\proj$ instead of $(\tilde v)_\proj$ from now on.

\begin{lemma}\label{lem:2}
\[
\displaystyle \tilde v_\proj^{(l)} =
\begin{cases}
\displaystyle \frac2d\bigg(1-\sum_{k=j+1}^{d-1}\lambda_k \bigg) & \displaystyle \text{ for } l=0,\dots,j,
\\[2ex]
\displaystyle \frac2d\bigg(1-\sum_{\substack{k=j+1\\k\not=l}}^{d-1}\lambda_k\bigg) + \big(1-\frac2d\big)\lambda_l & \displaystyle \text{ for } l=j+1,\dots,d-1,
\end{cases}
\]
where $\tilde v =p+\lambda_{j+1} (e_{j+1}-p) +\dots +\lambda_{d-1} (e_{d-1}-p) \in A''$.
\end{lemma}

\begin{proof}
From the proofs of Fact 1 and 3 it follows that\\[2ex]
\centerline{$
\tilde{v}_{\proj} \in A( \{ e_{j+1}, \ldots, e_{d-1}, \frac{2}{d} \1_d \} ) \cap \mathcal{S}_d.
$}~\\[2ex]
The smallest affine space, which contains this set is given by\\[2ex]
\centerline{$
A=\left\{ a \in \R^d \mid a = \frac{2}{d} \1_d + \delta_{j+1} (e_{j+1} - \frac{2}{d} \1_{d})
+ \ldots + \delta_{d-1} (e_{d-1} - \frac{2}{d} \1_d) \right\}.
$}~\\[2ex]
To find $\tilde{v}_{\proj}$ the projection for $v$ on $\mathcal{S}_d$, we minimize the distance of $v$ to $A$.

For $a \in A$
\begin{eqnarray} \label{eq:3.1}
 d(\tilde{v},a)^2
& = & \sum_{l=0}^{j} \Bigg( \tilde{v}^{(l)} - \frac{2}{d}
+ \delta_{j+1} \frac{2}{d}
+ \ldots + \delta_{d-1} \frac{2}{d} \Bigg)^2\\
& & + \sum_{l=j+1}^{d-1} \Bigg( \tilde{v}^{(l)} - \frac{2}{d} - \delta_l \left( 1 - \frac{2}{d} \right)
+ \sum_{\substack{k=j+1 \\ k \not= l} }^{d-1} \delta_k \frac{2}{d} \Bigg)^2 \nonumber.
\end{eqnarray}
 Calculating partial derivatives with respect to $\delta_{i}$, $i= j+1, \ldots , d-1$, yields

\begin{center}
$ \begin{array}{rcl}
\displaystyle \frac{\partial d (\tilde{v}, a)^2}{\partial \delta_{i}}
 &=& \displaystyle \sum_{l=0}^{j} 2 \Bigg( \tilde{v}^{(l)} - \frac{2}{d}
+ \delta_{j+1} \frac{2}{d} + \ldots + \delta_{d-1} \frac{2}{d} \Bigg)\frac{2}{d} \\
 & &\displaystyle + \sum_{l=j+1}^{d-1} 2 \Bigg( \tilde{v}^{(l)} -\frac{2}{d} - \delta_l \left( 1- \frac{2}{d} \right)
 + \sum_{\substack{k=j+1 \\k \not= l }}^{d-1} \delta_{k} \frac{2}{d} \Bigg) \alpha\\
 &=& \displaystyle 2 \Bigg( \frac{2}{d} \sum_{\substack{l=0 \\ l \not= i }}^{d-1} \tilde{v}^{(l)}
 - \left( 1 - \frac{2}{d} \right) \tilde{v}^{(i)} - \frac{2}{d} + \delta_{i} \Bigg) \ , \\[12mm]
\multicolumn{3}{l}{\mbox{where } \alpha = \frac{2}{d} \mbox{ for } l \not= i, \ \alpha = -(1-\frac{2}{d}) \mbox{ for } l = i, \mbox{ and thus} } \\[2ex]
\multicolumn{3}{l}{\displaystyle \frac{\partial d (\tilde{v}, a)^2}{\partial \delta_{i}} =0
\ \Leftrightarrow \ \delta_i = \frac{2}{d} \Bigg( 1 - \sum_{\substack{l=0 \\ l \not= i }}^{d-1} \tilde{v}^{(l)} \Bigg)  + \left( 1 - \frac{2}{d} \right) \tilde{v}^{(i)} .}
 \end{array}$
\end{center}

The determinant of the Hessian is positive which shows that a minimum occurs. According to the statement of Lemma \ref{lem:2} the components of $\tilde{v}$ has the following representation
\begin{eqnarray} \label{eq:3.2}
 \displaystyle \tilde{v}^{(l)} =
\begin{cases}
  \displaystyle p^{(l)} \left( 1 - \sum_{k = j+1}^{d-1} \lambda _k \right) \quad & \displaystyle \mbox{for } l = 0, \ldots , j, \\
  \displaystyle \lambda_l & \displaystyle \mbox{for } l= j+1, \ldots, d-1,
\end{cases}
\end{eqnarray}
where one should note that $p^{(j+1)} = \ldots = p^{(d-1)} = 0$.

Plugging in the equation of $\delta_i$, $i=j+1, \ldots, d-1$, and noting that
$\sum_{l=0}^{j} p^{(l)} = 1$ leads to
$$
\delta_i = \frac{2}{d} \left( 1 - \sum_{l=0}^{j} p ^{(l)} \left( 1 - \sum_{k=j+1}^{d-1} \lambda_k \right)
- \sum_{\substack{l=j+1\\ l \not= i}}^{d-1} \lambda_k \right) + \left( 1 - \frac{2}{d} \right) \lambda_i
$$
and finally to $\delta_i = \lambda_i$. Plugging this in equation (\ref{eq:3.1}) leads to the statement of the Lemma.
\end{proof}

\begin{lemma} \label{lem:3}
It holds:
\begin{enumerate}[1)]
 \item ~\vspace*{-7mm}

\refstepcounter{equation}%
\label{eq:3.3}
$ \displaystyle (\tilde{v} - \tilde{v}_{\proj} )^{(l)} =
\begin{cases}
  \displaystyle \left( p^{(l)} -\frac{2}{d} \right) \left( 1 - \sum_{k = j+1}^{d-1} \lambda _k \right) \quad & \displaystyle \mbox{for } l = 0, \ldots , j , \\
  \displaystyle - \frac{2}{d} \left( 1 - \sum_{k = j+1}^{d-1} \lambda _k \right) \quad & \displaystyle \mbox{for } l= j+1, \ldots, d-1.
\end{cases} \hfill (\arabic{section}.\arabic{equation}) $
\item The smallest affine subspace which contains $\mathcal{R} (p)$ can be
expressed as $x+ U$ where one can choose $x = \tilde{v}_{\proj}$ and
$$
\begin{array}{ll}
\displaystyle e_i + p^{(i)} \1_d - \tilde{v}_{\proj} &\displaystyle \quad \mbox{for } i=0, \ldots , j \\
\displaystyle e_i - \tilde{v}_{\proj} & \displaystyle \quad \mbox{for } i = j+1, \ldots, d-1
\end{array}
$$
as linear generating system of $U$.
\end{enumerate}
\end{lemma}

\begin{proof}
Statement 1) is a direct consequence of Lemma \ref{lem:2} and (\ref{eq:3.1}).
Statement 2) follows from the fact that $\tilde{v}_{\proj} = v_{\proj} \in \mathcal{R} (p(v))$
and that $\mathcal{R}(p) = \conv (e_i+ p^{(i)} \1_d \mid i=1,\ldots, d-1)$ where
$p^{(j+1)} = \ldots = p^{(d-1)} = 0$.

\begin{lemma} \label{lem:4}
 It holds
\[
\tilde{v} - \tilde{v}_{\proj} \bot e_i + p^{(i)} \1_d - \tilde{v}_{\proj} \quad \mbox{for } i = 0, \ldots, j.
\]
\end{lemma}

\begin{proof}
Lemma \ref{lem:2} implies
%
\begin{eqnarray}\label{eq:3.4}
\lefteqn{\left( e_i + p^{(i)} \1_d - \tilde{v}_{\proj} \right)^{(l)} } \qquad \nonumber\\[1ex]
&=&
\begin{cases}
p^{(i)} - \frac{2}{d} \left( 1 - \sum\limits_{k=j+1}^{d-1} \lambda_k \right) & l=0, \ldots, j; \ l\not= i,
\\
1+ p^{(i)} - \frac{2}{d} \left( 1- \sum\limits_{k=j+1}^{d-1} \lambda_k \right)  & l = i, 
\\
p^{(i)} - \frac{2}{d} \left( 1 - \sum\limits_{k=j+1, k \not= l}^{d-1} \lambda_k \right) - \left( 1- \frac{2}{d} \right) \lambda_l  & l= j+1 , \ldots, d-1.
\end{cases}
\end{eqnarray}
From (\ref{eq:3.3}) and (\ref{eq:3.4}) it follows
\begin{eqnarray*}
\lefteqn{\langle \tilde{v} - \tilde{v}_{\proj}, e_i + p^{(i)} \1_d - \tilde{v}_{\proj} \rangle}  \quad\\
&=& \sum_{\substack{l=0\\l \not= i}}^{j} \left( p^{(l)}
- \frac{2}{d} \right) \left( 1 - \sum\limits_{k=j+1}^{d-1} \lambda_k \right) \left( p^{(i)}
- \frac{2}{d} \left( 1- \sum\limits_{k=j+1}^{d-1} \lambda_k \right)  \right)
\\
&&{}+\left( p^{(i)}-\frac2d\right) \left(1-\sum\limits_{k=j+1}^{d-1}\lambda_k\right) \left(1+p^{(i)}-\frac2d\left(1-\sum\limits_{k=j+1}^{d-1}\lambda_k\right)\right)
\\
&&{}-\sum\limits_{l=j+1}^{d-1}\frac2d\left(1-\sum\limits_{k=j+1}^{d-1}\lambda_k\right)\left(p^{(i)} -\frac2d\left(1-\sum\limits_{\substack{k=j+1}{k\not= l}}^{d-1} \lambda_k\right) -\left(1-\frac2d\right)\lambda_l\right)
\\
&=& \left( 1 - \sum\limits_{k=j+1}^{d-1} \lambda_k \right) \cdot
 \Bigg[ \left( p^{(i)} - \frac{2}{d} \right)
+ \sum\limits_{i=0}^{j} p^{(l)}\left( p^{(i)} - \frac2d \left( 1- \sum\limits_{k=j+1}^{d-1} \lambda_k \right) \right)
\\
&& {} - \sum\limits_{l=0}^{i}\frac2d \left( p^{(i)} - \frac2d \left( 1 -  \sum\limits_{k=j+1}^{d-1} \lambda_k \right) \right)
- \sum\limits_{l=j+1}^{d-1} \frac{2}{d} \left( p^{(i)} - \frac{2}{d} \left( 1- \sum\limits_{k=j+1}^{d-1} \lambda_k \right) \right)
\\
&& {} + \sum\limits_{k=j+1}^{d-1} \frac2d \frac2d \lambda_k + \sum\limits_{l=j+1}^{d-1} \frac2d  \left( 1 -\frac{2}{d} \right) \lambda_l \Bigg]
\\
&=& \left( 1 - \sum\limits_{k=j+1}^{d-1} \lambda_k \right) \\
&&  \cdot \Bigg[ p^{(i)} - \frac2d + p^{(i)} - \frac2d \left( 1- \sum\limits_{k=j+1}^{d-1} \lambda_k \right)
- d \frac2d \left( p^{(i)} - \frac2d \left( 1- \sum\limits_{k=j+1}^{d-1} \lambda_k \right) \right)
\\
&& \quad {} + \sum\limits_{k=j+1}^{d-1} \frac2d \frac2d \lambda_k + \sum\limits_{l=j+1}^{d-1} \frac2d \lambda_l - \sum\limits_{l=j+1}^{d-1} \frac2d \frac2d \lambda_l \Bigg] = 0.\\[-9ex]
\end{eqnarray*}
\end{proof}

\medskip

\begin{lemma} \label{lem:5}
It holds
\[
\tilde{v} - \tilde{v}_{\proj}  \quad \bot \quad  e_i + \tilde{v}_{\proj} \qquad \mbox{for } i = j+1, \ldots, d-1.
\]
\end{lemma}

\newpage

\begin{proof}
By Lemma \ref{lem:2} one gets
%
\begin{eqnarray} \label{eq:3.5}
\lefteqn{\left( e_i- \tilde{v}_{\proj} \right)^{(l)}} \quad
\nonumber\\
&=& \begin{cases}
 - \frac{2}{d} \Bigg( 1 - \sum\limits_{\substack{k=j+1\\\phantom{k}}}^{d-1} \lambda_k \Bigg) & l=0, \ldots, j,
\\
 - \frac{2}{d} \Bigg( 1- \sum\limits_{\substack{k=j+1 \\ k \not= l}}^{d-1} \lambda_k \Bigg)
- \left( 1- \frac{2}{d} \right) \lambda_l & l= j+1 , \ldots, d-1;\ l\not= i,
\\
1 - \frac{2}{d} \Bigg( 1 - \sum\limits_{\substack{k=j+1 \\ k \not= i}}^{d-1} \lambda_k \Bigg) - \left( 1- \frac{2}{d} \right) \lambda_i & l = i.
 \end{cases}
\end{eqnarray}
From (\ref{eq:3.3}) and (\ref{eq:3.5}) it follows
\begin{eqnarray*}
\lefteqn{\langle \tilde{v} - \tilde{v}_{\proj}, e_i - \tilde{v}_{\proj} \rangle}\quad
\\
&=& \sum_{l=0 }^{j} \left( p^{(l)} - \frac{2}{d} \right) \left( 1 - \sum_{k=j+1}^{d-1} \lambda_k \right) \left(- \frac{2}{d} \left( 1- \sum_{k=j+1}^{d-1} \lambda_k \right)  \right)
\\
&& {} + \sum_{\substack{l=j+1\\l \not= i}}^{d-1} -\frac2d\ \left(1-\sum_{k=j+1}^{d-1}\lambda_k \right)
\left(-\frac2d \left(1-\sum_{\substack{k=j+1\\ k \not= l}}^{d-1}\lambda_k \right) - \left( 1- \frac2d\right)\lambda_l \right)
\\
&& {}- \frac2d \left( 1 - \sum_{k=j+1}^{d-1}\lambda_k\right) \left( 1 - \frac2d \left(1-\sum_{\substack{k=j+1\\ k\not= i}}^{d-1} \lambda_k \right) - \left(1-\frac2d\right)\lambda_i\right)
\\
&=& \left( 1 - \sum_{k=j+1}^{d-1} \lambda_k \right)  \\
&& \cdot \Bigg[ \sum_{l=0}^{i} p^{(l)} \left(- \frac2d \left( 1- \sum_{k=j+1}^{d-1} \lambda_k \right) \right) + \sum_{l=0}^{j} \frac2d \frac2d \left( 1 -  \sum_{k=j+1}^{d-1} \lambda_k \right)
 \\
&& \quad {} + \sum_{l=j+1}^{d-1} \frac2d \frac2d \left( 1- \sum_{k=j+1}^{d-1} \lambda_k \right) + \sum_{k=j+1}^{d-1} \frac2d \frac2d \lambda_k + \sum_{l=j+1}^{d-1} \frac{2}{d} \left( 1 -\frac{2}{d} \right) \lambda_l - \frac{2}{d} \Bigg]
\\
&=& \left( 1 - \sum_{k=j+1}^{d-1} \lambda_k \right) \\
&&  \cdot\Bigg[ - \frac2d \left( 1- \sum_{k=j+1}^{d-1} \lambda_k \right)
+ d \frac2d \frac2d \left( 1- \sum_{k=j+1}^{d-1} \lambda_k \right)
+ \sum_{k=j+1}^{d-1} \frac2d \frac2d \lambda_k
\\
&& \quad {}+ \sum_{l=j+1}^{d-1} \frac2d \lambda_l - \sum_{l=j+1}^{d-1} \frac2d \frac2d \lambda_l - \frac2d \Bigg]
= 0.\\[-8ex]
\end{eqnarray*}
\end{proof}

Finally we can state the proof of Fact 2: By Lemma \ref{lem:3}, \ref{lem:4}, and \ref{lem:5} one has
$\tilde{v} - \tilde{v}_{\proj} \ \bot \ A (\mathcal{R}(p))$.
By Fact 1 it follows that
$v - v_{\proj} \ \bot \ A (\mathcal{R}(p))$.
\end{proof}


\noindent%
\textbf{Acknowledgements.}

\end{document}